\documentclass[leqno,draft]{article}

\newtheorem{theorem}{Theorem}
\newtheorem{lemma}[theorem]{Lemma}
\newtheorem{proposition}[theorem]{Proposition}
\newtheorem{definition}[theorem]{Definition}
\newtheorem{corollary}[theorem]{Corollary}

\newcommand{\begintheorem}{\addtocounter{equation}{1}\begin{theorem}}
\newcommand{\beginlemma}{\addtocounter{equation}{1}\begin{lemma}}
\newcommand{\beginproposition}{\addtocounter{equation}{1}\begin{proposition}}
\newcommand{\begindefinition}{\addtocounter{equation}{1}\begin{definition}}
\newcommand{\begincorollary}{\addtocounter{equation}{1}\begin{corollary}}

\begin{document}

\title{$p$-Adic Heisenberg Cantor sets, 3}

\author{Stephen Semmes \\
        Rice University}

\date{}

\maketitle

\begin{abstract}
Here we look at some related constructions of solenoids, and mappings
associated to them.
\end{abstract}

\tableofcontents

\section{The Heisenberg groups $H_n({\bf R})$}
\label{heisenberg groups H_n(R)}
\setcounter{equation}{0}

        Let $n$ be a positive integer, let ${\bf R}$ be the real line,
and let $H_n({\bf R})$ be defined initially as a set as ${\bf R}^n
\times {\bf R}^n \times {\bf R}$.  More precisely, this is a smooth
manifold in the usual way, and a nice topological space in particular.
Elements of $H_n({\bf R})$ may typically be denoted as $(x, y, t)$,
with $x, y \in {\bf R}^n$ and $t \in {\bf R}$, so that $x$, $y$ have
components $x_j, y_j \in {\bf R}$ for $j = 1, \ldots, n$.

        If $(x, y, t), (x', y', t') \in H_n({\bf R})$, then put
\begin{equation}
\label{(x, y, t) diamond (x', y', t') = (x + x', y + y', t + t' + x cdot y')}
 (x, y, t) \diamond (x', y', t') = (x + x', y + y', t + t' + x \cdot y'),
\end{equation}
where $x \cdot y'$ refers to the usual dot product on ${\bf R}^n$, so that
\begin{equation}
        x \cdot y' = \sum_{j = 1}^n x_j \, y_j'.
\end{equation}
If $(x'', y'', t'') \in H_n({\bf R})$ too, then
\begin{eqnarray}
\label{((x, y, t) diamond (x', y', t')) diamond (x'', y'', t'') = ...}
\lefteqn{((x, y, t) \diamond (x', y', t')) \diamond (x'', y'', t'')} \\
 & = & (x + x', y + y', t + t' + x \cdot y') \diamond (x'', y'', t'')
                                                                   \nonumber \\
 & = & (x + x' + x'', y + y' + y'', t + t' + t'' + x \cdot y'
                                     + x \cdot y'' + x' \cdot y''). \nonumber
\end{eqnarray}
Similarly,
\begin{eqnarray}
\label{(x, y, t) diamond ((x', y', t') diamond (x'', y'', t'')) = ...}
\lefteqn{(x, y, t) \diamond ((x', y', t') \diamond (x'', y'', t''))} \\
 & = & (x, y, t) \diamond (x' + x'', y' + y'', t' + t'' + x' \cdot y'')
                                                                   \nonumber \\
 & = & (x + x' + x'', y + y' + y'', t + t' + t'' + x \cdot y' + x \cdot y''
                                                  + x' \cdot y''). \nonumber
\end{eqnarray}
Thus
\begin{equation}
\label{associativity}
        ((x, y, t) \diamond (x', y', t')) \diamond (x'', y'', t'')
            = (x, y, t) \diamond ((x', y', t') \diamond (x'', y'', t''),
\end{equation}
which shows that $\diamond$ is associative on $H_n({\bf R})$.  Of
course, $\diamond$ is not commutative on $H_n({\bf R})$, because of
the $x \cdot y'$ term in (\ref{(x, y, t) diamond (x', y', t') = (x +
x', y + y', t + t' + x cdot y')}).

        Note that
\begin{equation}
\label{(0, 0, 0) diamond (x, y, t) = (x, y, t) diamond (0, 0, 0) = (x, y, t)}
        (0, 0, 0) \diamond (x, y, t) = (x, y, t) \diamond (0, 0, 0) = (x, y, t)
\end{equation}
for every $(x, y, t) \in H_n({\bf R})$, so that $(0, 0, 0)$ is the
identity element of $H_n({\bf R})$ with respect to $\diamond$.  Put
\begin{equation}
\label{(x, y, t)^{-1} = (-x, -y, -t + x cdot y)}
        (x, y, t)^{-1} = (-x, -y, -t + x \cdot y),
\end{equation}
for each $(x, y, t) \in H_n({\bf R})$, and observe that
\begin{equation}
\label{(x, y, t) diamond (x, y, t)^{-1} = ... = (0, 0, 0)}
        (x, y, t) \diamond (x, y, t)^{-1}
           = (x, y, t) \diamond (-x, -y, -t + x \cdot y) = (0, 0, 0),
\end{equation}
and similarly that
\begin{equation}
\label{(x, y, t)^{-1} diamond (x, y, t) = ... = (0, 0, 0)}
        (x, y, t)^{-1} \diamond (x, y, t)
           = (-x, -y, -t + x \cdot y) \diamond (x, y, t) = (0, 0, 0).
\end{equation}
Hence $H_n({\bf R})$ is a group with respect to $\diamond$, and $(x,
y, t)^{-1}$ is the inverse of $(x, y, t)$ in this group.  More
precisely, $H_n({\bf R})$ is a Lie group, since the group operations
are given by smooth mappings.

        Consider the mapping $\pi : H_n({\bf R}) \to {\bf R}^n \times
{\bf R}^n$ defined by
\begin{equation}
        \pi((x, y, t)) = (x, y).
\end{equation}
This is a homomorphism from $H_n({\bf R})$ as a group with respect to
$\diamond$ onto ${\bf R}^n \times {\bf R}^n$ as a commutative group
with respect to coordinatewise addition, because
\begin{equation}
        \pi((x, y, t) \diamond (x', y', t') = (x + x', y + y')
\end{equation}
and
\begin{equation}
        \pi((x, y, t)^{-1}) = (-x, -y)
\end{equation}
for every $(x, y, t), (x', y', t') \in H_n({\bf R})$.  The kernel of
this homomorphism is the subgroup of $H_n({\bf R})$ consisting of
elements of the form $(0, 0, t)$, $t \in {\bf R}$, which is in fact
the center of $H_n({\bf R})$.  It is easy to see that the commutator
of any two elements of $H_n({\bf R})$ is in the center, which implies
that $H_n({\bf R})$ is a nilpotent group.  Equivalently, the
commutator of any two elements of $H_n({\bf R})$ is in the kernel of
$\pi$, since $\pi$ is a homomorphism from $H_n({\bf R})$ onto a
commutative group.

\section{Some subgroups of $H_n({\bf R})$}
\label{some subgroups of H_n(R)}
\setcounter{equation}{0}

        Let $A$ be a subring of the ring ${\bf R}$ of real numbers.
Thus $A$ is a subgroup of ${\bf R}$ as a commutative group with
respect to addition, and the product of any two elements of $A$ is
also an element of $A$.  Note that $1$ is not required to be an
element of $A$.  Under these conditions, it is easy to see that
\begin{equation}
        H_n(A) = A^n \times A^n \times A
\end{equation}
is a subgroup of $H_n({\bf R})$ with respect to $\diamond$.  However,
$H_n(A)$ is not a normal subgroup of $H_n({\bf R})$, except in the
trivial cases where $A = {\bf R}$ or $A = \{0\}$.  Indeed, if $(x, y,
t), (x', y', t') \in H_n({\bf R})$, then
\begin{eqnarray}
\label{((x', y', t') diamond (x, y, t)) diamond (x', y', t')^{-1} = ...}
\lefteqn{((x', y', t') \diamond (x, y, t)) \diamond (x', y', t')^{-1}} \\
 & = & (x' + x, y' + y, t' + t + x' \cdot y) \diamond
                              (- x', - y', - t' + x' \cdot y') \nonumber \\
 & = & (x, y, t + x' \cdot y - x' \cdot y' - x \cdot y' + x' \cdot y')
                                                                \nonumber \\
 & = & (x, y, t + x' \cdot y - x \cdot y'). \nonumber
\end{eqnarray}
If $A \ne \{0\}$, then one can choose $(x, y, t) \in H_n(A)$, $(x',
y', t') \in H_n({\bf R})$ so that
\begin{equation}
\label{t + x' cdot y - x cdot y'}
        t + x' \cdot y - x \cdot y'
\end{equation}
is any real number, and hence not an element of $A$ unless $A = {\bf
R}$.  If $B$ is a subring of ${\bf R}$ which is an ideal in $A$, then
$H_n(B)$ is a normal subgroup of $H_n(A)$.

        If $A$ is a subgroup of ${\bf R}$ as a commutative group with
respect to addition, then the quotient ${\bf R} / A$ can be defined as
an abelian group in the usual way.  Similarly, $A^n \times A^n$ may be
considered as an abelian subgroup of ${\bf R}^n \times {\bf R}^n$ with
respect to coordinatewise addition, and the quotient group
\begin{equation}
\label{({bf R}^n times {bf R}^n) / (A^n times A^n)}
        ({\bf R}^n \times {\bf R}^n) / (A^n \times A^n)
\end{equation}
is isomorphic to
\begin{equation}
\label{({bf R} / A)^n times ({bf R} / A)^n}
        ({\bf R} / A)^n \times ({\bf R} / A)^n.
\end{equation}
As in the previous section, there is a natural homomorphism from
$H_n({\bf R})$ onto ${\bf R}^n \times {\bf R}^n$, which can then be
composed with the quotient mapping from ${\bf R}^n \times {\bf R}^n$
onto (\ref{({bf R}^n times {bf R}^n) / (A^n times A^n)}) to get a
homomorphism from $H_n({\bf R})$ onto (\ref{({bf R}^n times {bf R}^n)
/ (A^n times A^n)}).  The kernel of this homomorphism from $H_n({\bf
R})$ onto (\ref{({bf R}^n times {bf R}^n) / (A^n times A^n)}) is equal to
\begin{equation}
\label{A^n times A^n times {bf R}}
        A^n \times A^n \times {\bf R},
\end{equation}
which is a normal subgroup of $H_n({\bf R})$.

        Let $A$ be a subring of ${\bf R}$ again, and observe that
$H_n(A)$ is a normal subgroup of (\ref{A^n times A^n times {bf R}})
with respect to $\diamond$.  The quotient of (\ref{A^n times A^n times
{bf R}}) by $H_n(A)$ is isomorphic to the abelian group ${\bf R} / A$
in the obvious way.  If $A \ne \{0\}$, then (\ref{A^n times A^n times
{bf R}}) is the ``normal closure'' of $H_n(A)$ in $H_n({\bf R})$,
which is the smallest normal subgroup of $H_n({\bf R})$ that contains
$H_n(A)$.  This follows from the fact that (\ref{t + x' cdot y - x
cdot y'}) can be any real number when $(x, y, t) \in H_n(A)$ and $(x',
y', t') \in H_n({\bf R})$, as before.  Even if $H_n(A)$ is not a
normal subgroup of $H_n({\bf R})$, we can still consider the quotient
space $H_n({\bf R}) / H_n(A)$ as a set, with the canonical quotient
mapping from $H_n({\bf R})$ onto $H_n({\bf R}) / H_n(A)$.
There is also a natural mapping from $H_n({\bf R}) / H_n(A)$ onto
\begin{equation}
\label{H_n({bf R}) / (A^n times A^n times {bf R}) cong ...}
        H_n({\bf R}) / (A^n \times A^n \times {\bf R})
         \cong ({\bf R}^n \times {\bf R}^n) / (A^n \times A^n),
\end{equation}
because $H_n(A)$ is a subgroup of (\ref{A^n times A^n times {bf R}}).
The quotient mapping from $H_n({\bf R})$ onto (\ref{H_n({bf R}) / (A^n
times A^n times {bf R}) cong ...})  is the same as the composition of
the quotient mapping from $H_n({\bf R})$ onto $H_n({\bf R}) / H_n(A)$
with this new mapping from $H_n({\bf R}) / H_n(A)$ onto (\ref{H_n({bf
R}) / (A^n times A^n times {bf R}) cong ...}).

        Although $H_n({\bf R}) / H_n(A)$ is not a group when $A \ne
{\bf R}, \{0\}$, there is still a natural action of $H_n({\bf R})$ on
$H_n({\bf R}) / H_n(A)$, corresponding to multiplication by elements
of $H_n({\bf R})$ on the left.  This is because one can multiply a
left coset of $H_n(A)$ in $H_n({\bf R})$ by an element of $H_n({\bf
R})$ on the left and get another left coset of $H_n(A)$ in $H_n({\bf
R})$.  The quotient mapping from $H_n({\bf R})$ onto $H_n({\bf R}) /
H_n(A)$ intertwines the obvious action of $H_n({\bf R})$ on itself by
left-multiplication with this action on $H_n({\bf R}) / H_n(A)$.  The
analogous action of $H_n({\bf R})$ on the quotient group (\ref{H_n({bf
R}) / (A^n times A^n times {bf R}) cong ...}) is the same as first
mapping $H_n({\bf R})$ onto (\ref{H_n({bf R}) / (A^n times A^n times
{bf R}) cong ...}) and then using the group operation on the quotient.
The mapping from $H_n({\bf R}) / H_n(A)$ onto (\ref{H_n({bf R}) / (A^n
times A^n times {bf R}) cong ...}) described in the preceding
paragraph intertwines these actions of $H_n({\bf R})$ on $H_n({\bf R})
/ H_n(A)$ and (\ref{H_n({bf R}) / (A^n times A^n times {bf R}) cong ...}).

\section{Some discrete subgroups}
\label{some discrete subgroups}
\setcounter{equation}{0}

        The ring ${\bf Z}$ of integers is a subring of ${\bf R}$, and
a subgroup of ${\bf R}$ as a commutative group with respect to
addition in particular.  Thus ${\bf R} / {\bf Z}$ may be considered as
a commutative group as well.  One can also consider ${\bf R} / {\bf
Z}$ as a topological space, and indeed as a $1$-dimensional smooth
manifold, such that the quotient mapping from ${\bf R}$ onto ${\bf R}
/ {\bf Z}$ is a local diffeomorphism.  The group operations on ${\bf
R} / {\bf Z}$ are smooth, so that ${\bf R} / {\bf Z}$ is actually a
Lie group.  It is well known that ${\bf R} / {\bf Z}$ is isomorphic as
a Lie group to the group ${\bf T}$ of complex numbers $z$ with $|z| =
1$ with respect to multiplication.

        Similarly, ${\bf Z}^n$ is a subgroup of ${\bf R}^n$ as a
commutative group with respect to coordinatewise addition, and the
quotient ${\bf R}^n / {\bf Z}^n$ is a commutative group isomorphic to
$({\bf R} / {\bf Z})^n \cong {\bf T}^n$.  One can also consider ${\bf
R}^n / {\bf Z}^n$ as an $n$-dimensional smooth manifold, for which the
quotient mapping from ${\bf R}^n$ onto ${\bf R}^n / {\bf Z}^n$ is a
local diffeomorphism, and which is diffeomorphic to the
$n$-dimensional torus ${\bf T}^n$.  The group operations on ${\bf R}^n
/ {\bf Z}^n$ are smooth, so that ${\bf R}^n / {\bf Z}^n$ is a Lie
group, which is isomorphic to ${\bf T}^n$.

        As in the previous section, $H_n({\bf Z})$ is a subgroup of
$H_n({\bf R})$, but not a normal subgroup.  The quotient space
$H_n({\bf R}) / H_n({\bf Z})$ can still be defined as a set, and it
may also be considered as a compact smooth manifold of dimension $2 n
+ 1$, for which the corresponding quotient mapping from $H_n({\bf R})$
onto $H_n({\bf R}) / H_n({\bf Z})$ is a local diffeomorphism.
Remember that $H_n({\bf R})$ is the same as ${\bf R}^n \times {\bf
R}^n \times {\bf R}$ as a smooth manifold.  Although $H_n({\bf R}) /
H_n({\bf Z})$ is not a group, there is a natural action of $H_n({\bf
R})$ on $H_n({\bf R}) / H_n({\bf Z})$ by multiplication on the left,
as before.  It is easy to see that this action corresponds to a smooth
mapping from the Cartesian product of $H_n({\bf R})$ and $H_n({\bf R})
/ H_n({\bf Z})$ into $H_n({\bf R}) / H_n({\bf Z})$, because of the
smoothness of the group operation on $H_n({\bf R})$.

        The quotient group (\ref{H_n({bf R}) / (A^n times A^n times
{bf R}) cong ...}) with $A = {\bf Z}$ is the same as
\begin{equation}
\label{({bf R}^n times {bf R}^n) / ({bf Z}^n times {bf Z}^n) cong ...}
        ({\bf R}^n \times {\bf R}^n) / ({\bf Z}^n \times {\bf Z}^n)
                                          \cong {\bf T}^n \times {\bf T}^n,
\end{equation}
which may be considered as a $2 n$-dimensional smooth manifold and a
Lie group.  The corresponding quotient mapping from $H_n({\bf R})$
onto ${\bf T}^n \times {\bf T}^n$ is a smooth submersion.  This
basically means that the quotient mapping looks locally like the
standard projection from ${\bf R}^n \times {\bf R}^n \times {\bf R}$
onto ${\bf R}^n \times {\bf R}^n$, which is clear from the
construction.  The natural mapping from $H_n({\bf R}) / H_n({\bf Z})$
onto (\ref{H_n({bf R}) / (A^n times A^n times {bf R}) cong ...}) with
$A = {\bf Z}$ discussed in the previous section can be identified with
a mapping from $H_n({\bf R}) / H_n({\bf Z})$ onto ${\bf T}^n \times
{\bf T}^n$, which is a smooth submersion as well, for the same reasons.

        It is sometimes helpful to look at
\begin{equation}
\label{{0} times {0} times {bf Z}}
        \{0\} \times \{0\} \times {\bf Z}
\end{equation}
as a discrete subgroup of the center of $H_n({\bf R})$, and the
corresponding quotient
\begin{equation}
\label{H_n({bf R}) / ({0} times {0} times {bf Z})}
        H_n({\bf R}) / (\{0\} \times \{0\} \times {\bf Z}).
\end{equation}
As usual, (\ref{H_n({bf R}) / ({0} times {0} times {bf Z})}) may be
considered as a $(2n + 1)$-dimensional smooth manifold, such that the
quotient mapping from $H_n({\bf R})$ onto (\ref{H_n({bf R}) / ({0}
times {0} times {bf Z})}) is a local diffeomorphism.  It is easy to
see that (\ref{H_n({bf R}) / ({0} times {0} times {bf Z})}) is
diffeomorphic to ${\bf R}^n \times {\bf R}^n \times {\bf T}$.  The
quotient (\ref{H_n({bf R}) / ({0} times {0} times {bf Z})}) is also a
group, because (\ref{{0} times {0} times {bf Z}}) is a normal subgroup
of $H_n({\bf R})$, and the quotient group is not commutative.  There
is a natural map from (\ref{H_n({bf R}) / ({0} times {0} times {bf
Z})}) onto $H_n({\bf R}) / H_n({\bf Z})$, since (\ref{{0} times {0}
times {bf Z}}) is a subgroup of $H_n({\bf Z})$, and this mapping is a
local diffeomorphism.

        Because (\ref{{0} times {0} times {bf Z}}) is a subgroup of
the center $\{0\} \times \{0\} \times {\bf R}$ of $H_n({\bf R})$,
there is a natural homomorphism from (\ref{H_n({bf R}) / ({0} times
{0} times {bf Z})}) onto
\begin{equation}
\label{H_n({bf R}) / ({0} times {0} times {bf R}) cong {bf R}^n times {bf R}^n}
        H_n({\bf R}) / (\{0\} \times \{0\} \times {\bf R})
                               \cong {\bf R}^n \times {\bf R}^n.
\end{equation}
This mapping is also a smooth submersion, which can be identified with
the standard projection from ${\bf R}^n \times {\bf R}^n \times {\bf
T}$ onto ${\bf R}^n \times {\bf R}^n$ as smooth manifolds.
The kernel of this homomorphism is the subgroup
\begin{equation}
\label{({0} times {0} times {bf R}) / ({0} times {0} times {bf Z})}
 (\{0\} \times \{0\} \times {\bf R}) / (\{0\} \times \{0\} \times {\bf Z})
\end{equation}
of (\ref{H_n({bf R}) / ({0} times {0} times {bf Z})}), which is
isomorphic to ${\bf R} / {\bf Z} \cong {\bf T}$ as a group.  If we
identify (\ref{H_n({bf R}) / ({0} times {0} times {bf Z})}) with ${\bf
R}^n \times {\bf R}^n \times {\bf T}$ as a manifold, then (\ref{({0}
times {0} times {bf R}) / ({0} times {0} times {bf Z})}) corresponds
to the submanifold $\{0\} \times \{0\} \times {\bf T}$.

        Of course, there is also a natural homomorphism from
(\ref{H_n({bf R}) / ({0} times {0} times {bf R}) cong {bf R}^n times
{bf R}^n}) onto (\ref{({bf R}^n times {bf R}^n) / ({bf Z}^n times {bf
Z}^n) cong ...}), which is a local diffeomorphism as well.  The
composition of this homomorphism with the one from (\ref{H_n({bf R}) /
({0} times {0} times {bf Z})}) onto (\ref{H_n({bf R}) / ({0} times {0}
times {bf R}) cong {bf R}^n times {bf R}^n}) is the same as the
natural homomorphism from (\ref{H_n({bf R}) / ({0} times {0} times {bf
Z})}) onto ${\bf T}^n \times {\bf T}^n$ that results from the
inclusion of (\ref{{0} times {0} times {bf Z}}) in the normal subgroup
\begin{equation}
\label{{bf Z}^n times {bf Z}^n times {bf R}}
        {\bf Z}^n \times {\bf Z}^n \times {\bf R}
\end{equation}
of $H_n({\bf R})$.  The kernel of this homomorphism from (\ref{H_n({bf
R}) / ({0} times {0} times {bf Z})}) onto ${\bf T}^n \times {\bf T}^n$
is the subgroup
\begin{equation}
\label{({bf Z}^n times {bf Z}^n times {bf R}) / ({0} times {0} times {bf Z})}
        ({\bf Z}^n \times {\bf Z}^n \times {\bf R})
                         / (\{0\} \times \{0\} \times {\bf Z})
\end{equation}
of (\ref{H_n({bf R}) / ({0} times {0} times {bf Z})}).  If we identify
(\ref{H_n({bf R}) / ({0} times {0} times {bf Z})}) with ${\bf R}^n
\times {\bf R}^n \times {\bf T}$ as a manifold, then (\ref{({bf Z}^n
times {bf Z}^n times {bf R}) / ({0} times {0} times {bf Z})})
corresponds to the subset
\begin{equation}
\label{{bf Z}^n times {bf Z}^n times {bf T}}
        {\bf Z}^n \times {\bf Z}^n \times {\bf T}
\end{equation}
of this manifold.  Note that (\ref{({bf Z}^n times {bf Z}^n times {bf
R}) / ({0} times {0} times {bf Z})}) is also isomorphic as a group
which is a subgroup of (\ref{H_n({bf R}) / ({0} times {0} times {bf
Z})}) to (\ref{{bf Z}^n times {bf Z}^n times {bf T}}) as a commutative
group, where the group operations are defined coordinatewise.  
This is because the $x \cdot y'$ and $x \cdot y$ terms in (\ref{(x, y,
t) diamond (x', y', t') = (x + x', y + y', t + t' + x cdot y')}) and
(\ref{(x, y, t)^{-1} = (-x, -y, -t + x cdot y)}) are integers when $x,
y, y' \in {\bf Z}^n$, and hence reduce to $0$ in ${\bf R} / {\bf Z}
\cong {\bf T}$.

        It is easy to check that (\ref{{0} times {0} times {bf Z}}) is
the center of $H_n({\bf Z})$, and that the quotient of $H_n({\bf Z})$
by (\ref{{0} times {0} times {bf Z}}) is isomorphic to ${\bf Z}^n
\times {\bf Z}^n$ as a commutative group with respect to
coordinatewise addition.  Because (\ref{{0} times {0} times {bf Z}})
is contained in $H_n({\bf Z})$, there is a natural map from
(\ref{H_n({bf R}) / ({0} times {0} times {bf Z})}) onto $H_n({\bf R})
/ H_n({\bf Z})$, which is a local diffeomorphism.  Equivalently,
one can identify $H_n({\bf R}) / H_n({\bf Z})$ with the quotient of
(\ref{H_n({bf R}) / ({0} times {0} times {bf Z})}) by the discrete subgroup
\begin{equation}
\label{H_n({bf Z}) / ({0} times {0} times {bf Z})}
        H_n({\bf Z}) / (\{0\} \times \{0\} \times {\bf Z}).
\end{equation}
The composition of this mapping from (\ref{H_n({bf R}) / ({0} times
{0} times {bf Z})}) onto $H_n({\bf R}) / H_n({\bf Z})$ with the usual
one from $H_n({\bf R}) / H_n({\bf Z})$ onto ${\bf T}^n \times {\bf
T}^n$ is the same as the natural homomorphism from (\ref{H_n({bf R}) /
({0} times {0} times {bf Z})}) onto ${\bf T}^n \times {\bf T}^n$
mentioned in the preceding paragraph.

        At any rate, $H_n({\bf R}) / H_n({\bf Z})$ is a circle bundle
over ${\bf T}^n \times {\bf T}^n$, with the mapping from $H_n({\bf R})
/ H_n({\bf Z})$ onto (\ref{H_n({bf R}) / (A^n times A^n times {bf R})
cong ...}) with $A = {\bf Z}$ as the projection onto the base.  The
fundamental group of $H_n({\bf R}) / H_n({\bf Z})$ is isomorphic to
$H_n({\bf Z})$, because $H_n({\bf R})$ is the same as ${\bf R}^n
\times {\bf R}^n \times {\bf R}$ as a manifold, which is
simply-connected.  By contrast, ${\bf T}^n \times {\bf T}^n \times
{\bf T}$ is the trivial circle bundle over ${\bf T}^n \times {\bf
T}^n$, using the obvious projection.  This can also be represented as
the quotient of ${\bf R}^n \times {\bf R}^n \times {\bf R}$ as a
commutative group with respect to coordinatewise addition by the discrete
subgroup ${\bf Z}^n \times {\bf Z}^n \times {\bf Z}$, and the fundamental
group is isomorphic to ${\bf Z}^n \times {\bf Z}^n \times {\bf Z}$.
Of course, $H_n({\bf Z})$ is not commutative, and hence is not isomorphic
as a group to ${\bf Z}^n \times {\bf Z}^n \times {\bf Z}$.

\section{Dilations}
\label{dilations}
\setcounter{equation}{0}

        Let $r$ be a real number, and let $\delta_r : H_n({\bf R}) \to
H_n({\bf R})$ be defined by
\begin{equation}
\label{delta_r((x, y, t)) = (r x, r y, r^2 2)}
        \delta_r((x, y, t)) = (r \, x, r \, y, r^2 \, 2),
\end{equation}
where $r \, x, r \, y \in {\bf R}^n$ are defined coordinatewise, as usual.
Thus
\begin{equation}
\label{delta_r((0, 0, 0)) = (0, 0, 0)}
        \delta_r((0, 0, 0)) = (0, 0, 0),
\end{equation}
and it is easy to check that
\begin{equation}
\label{delta_r((x, y, t) diamond (x', y', t')) = ...}
        \delta_r((x, y, t) \diamond (x', y', t')) =
                       \delta_r((x, y, t)) \diamond \delta_r((x', y', t'))
\end{equation}
and
\begin{equation}
\label{delta_r((x, y, t)^{-1}) = (delta_r((x, y, t)))^{-1}}
        \delta_r((x, y, t)^{-1}) = (\delta_r((x, y, t)))^{-1}
\end{equation}
for every $(x, y, t), (x', y', t') \in H_n({\bf R})$, so that $\delta_r$
is a homomorphism from $H_n({\bf R})$ into itself.  Note that $\delta_r$
is the identity mapping on $H_n({\bf R})$ when $r = 1$, and that
\begin{equation}
\label{delta_r circ delta_{r'} = delta_{r r'}}
        \delta_r \circ \delta_{r'} = \delta_{r \, r'}
\end{equation}
for every $r, r' \in {\bf R}$.  If $r \ne 0$, then $\delta_r$ is an
automorphism of $H_n({\bf R})$, whose inverse is equal to
$\delta_{r^{-1}}$.

        Now let $r \ge 2$ be an integer, so that $\delta_r(H_n({\bf
Z}))$ is a subgroup of $H_n({\bf Z})$.  More precisely,
$\delta_r(H_n({\bf Z}))$ is a subgroup of $H_n(r \, {\bf Z})$, and in fact
\begin{equation}
\label{H_n(r^2 {bf Z}) subseteq delta_r(H_n({bf Z})) subseteq H_n(r {bf Z})}
        H_n(r^2 \, {\bf Z}) \subseteq \delta_r(H_n({\bf Z}))
                             \subseteq H_n(r \, {\bf Z}).
\end{equation}
Observe that $H_n(k \, {\bf Z})$ is a normal subgroup of $H_n({\bf
Z})$ for each positive integer $k$, but that $\delta_r(H_n({\bf Z}))$
is not a normal subgroup of $H_n({\bf Z})$.  This follows easily from
(\ref{((x', y', t') diamond (x, y, t)) diamond (x', y', t')^{-1} =
...}), and we also get that $\delta_r(H_n({\bf Z}))$ is a normal
subgroup of $H_n(r \, {\bf Z})$.  Similarly,
\begin{equation}
\label{(r {bf Z})^n times (r {bf Z})^n times {bf Z}}
        (r \, {\bf Z})^n \times (r \, {\bf Z})^n \times {\bf Z}
\end{equation}
is a normal subgroup of $H_n({\bf Z})$, and $\delta_r(H_n({\bf Z}))$
is a normal subgroup of (\ref{(r {bf Z})^n times (r {bf Z})^n times {bf Z}}).

        It is easy to see that
\begin{equation}
\label{{0} times {0} times (k {bf Z})}
        \{0\} \times \{0\} \times (k \, {\bf Z})
\end{equation}
is the center of $H_n(k \, {\bf Z})$ for each positive integer $k$, and that
the commutator of any two elements of $H_n(k \, {\bf Z})$ is contained in
\begin{equation}
\label{{0} times {0} times (k^2 {bf Z})}
        \{0\} \times \{0\} \times (k^2 \, {\bf Z}).
\end{equation}
More precisely, (\ref{{0} times {0} times (k^2 {bf Z})}) is the
commutator subgroup of $H_n(k \, {\bf Z})$, which is the subgroup
generated by commutators of elements of $H_n(k \, {\bf Z})$.  Thus the
quotient of the center of $H_n(k \, {\bf Z})$ by the commutator
subgroup of $H_n(k \, {\bf Z})$ is isomorphic to
\begin{equation}
\label{k {bf Z} / k^2 {bf Z} cong {bf Z} / k {bf Z}}
        k \, {\bf Z} / k^2 \, {\bf Z} \cong {\bf Z} / k \, {\bf Z}.
\end{equation}
This shows that $H_n(k \, {\bf Z})$ is not isomorphic to $H_n(k' \,
{\bf Z})$ when $k$, $k'$ are distinct positive integers, and in
particular that $H_n(k \, {\bf Z})$ is not isomorphic to $H_n({\bf
Z})$ when $k \ge 2$.  Of course, $\delta_r(H_n({\bf Z}))$ is
isomorphic to $H_n({\bf Z})$, by construction.

        As in the preceding section, the quotient space $H_n({\bf R})
/ H_n(k \, {\bf Z})$ may be considered as a compact smooth manifold of
dimension $2 n + 1$ for each positive integer $k$, in such a way that
the corresponding quotient mapping from $H_n({\bf R})$ onto $H_n({\bf
R}) / H_n(k \, {\bf Z})$ is a local diffeomorphism.  The fundamental
group of $H_n({\bf R}) / H_n(k \, {\bf Z})$ is isomorphic to $H_n(k \,
{\bf Z})$, because $H_n({\bf R})$ is simply-connected.  In particular,
the fundamental group of $H_n({\bf R}) / H_n(k \, {\bf Z})$ is not
isomorphic to the fundamental group of $H_n({\bf R}) / H_n({\bf Z})$
when $k \ge 2$, by the remarks in the previous paragraph.  In the same
way, $H_n({\bf R}) / \delta_r(H_n({\bf Z}))$ may be considered as a
compact smooth manifold of dimension $2 n + 1$, for which the
corresponding quotient mapping from $H_n({\bf R})$ onto $H_n({\bf R})
/ \delta_r(H_n({\bf Z}))$ is a local diffeomorphism.  It is easy to
see that $H_n({\bf R}) / \delta_r(H_n({\bf Z}))$ is diffeomorphic to
$H_n({\bf R}) / H_n({\bf Z})$, using the fact that $\delta_r$ is a
diffeomorphism from $H_n({\bf R})$ onto itself when $r \ne 0$.

\section{Invariant metrics}
\label{invariant metrics}
\setcounter{equation}{0}

        As before, $H_n({\bf R})$ is the same as ${\bf R}^n \times
{\bf R}^n \times {\bf R}$ as a smooth manifold, and hence the tangent
space to $H_n({\bf R})$ at $(0, 0, 0)$ is isomorphic to ${\bf R}^n
\times {\bf R}^n \times {\bf R}$ as a vector space over the real
numbers.  If one chooses an inner product on this vector space, then
there is a unique smooth Riemannian metric on $H_n({\bf R})$ which is
invariant under right translations and which agrees with the given
inner product at $(0, 0, 0)$.  Because of invariance under right
translations, this leads to a Riemannian metric on $H_n({\bf R}) /
H_n({\bf Z})$ such that the quotient mapping from $H_n({\bf R})$ onto
$H_n({\bf R}) / H_n({\bf Z})$ preserves the Riemannian metric at each
point.  Similar remarks apply to other quotients of $H_n({\bf R})$ by
discrete subgroups on the right, such as $H_n({\bf R}) / H_n(k \, {\bf
Z})$ for each positive integer $k$, and $H_n({\bf R}) /
\delta_r(H_n({\bf Z}))$ for each $r \ne 0$.

        Let us continue to identify the tangent space of $H_n({\bf
R})$ at $(0, 0, 0)$ with ${\bf R}^n \times {\bf R}^n \times \{\bf
Z\}$, and consider the subspace
\begin{equation}
\label{{bf R}^n times {bf R}^n times {0}}
        {\bf R}^n \times {\bf R}^n \times \{0\}.
\end{equation}
This leads to a unique smooth distribution of hyperplanes in the
tangent spaces at arbitrary points in $H_n({\bf R})$, which is
invariant under right translations on $H_n({\bf R})$, and which is
equal to (\ref{{bf R}^n times {bf R}^n times {0}}) at $(0, 0, 0)$.
These hyperplanes are called the \emph{horizontal} subspaces of the
corresponding tangent spaces of $H_n({\bf R})$, and one can check that
they are also invariant under the dilations $\delta_r$ for each $r \in
{\bf R}$.  More precisely, the differential of a dilation $\delta_r$
at $(0, 0, 0)$ maps the tangent space to $H_n({\bf R})$ at $(0, 0, 0)$
to itself, and it is easy to see that the differential of $\delta_r$
at $(0, 0, 0)$ also maps the horizontal subspace (\ref{{bf R}^n times
{bf R}^n times {0}}) at $(0, 0, 0)$ into itself.  In order to deal
with other points in $H_n({\bf R})$, one can use the fact that the
dilations $\delta_r$ define homomorphisms on $H_n({\bf R})$.

        A continuously-differentiable curve in $H_n({\bf R})$ is said
to be \emph{horizontal} if the tangent vector of the curve at any
point along the curve is contained in the horizontal subspace of the
tangent space of $H_n({\bf R})$ at that point.  It is well known that
any pair of points in $H_n({\bf R})$ can be connected by a horizontal
curve, so that the distribution of horizontal subspaces of the tangent
spaces of $H_n({\bf R})$ is completely non-integrable.  If $H_n({\bf
R})$ is equipped with a Riemannian metric, then the length of a $C^1$
curve in $H_n({\bf R})$ can be defined in the usual way.  This leads
to the associated Riemannian distance between two points, which is the
infimum of the lengths of the curves connecting the two points.
Similarly, the sub-Riemannian distance between two points in $H_n({\bf
R})$ associated to this distribution of horizontal subspaces of the
tangent spaces is defined to be the infimum of the lengths of the
horizontal curves connecting the two points.  Thus the sub-Riemannian
distance is automatically greater than or equal to the Riemannian
distance corresponding to the same Riemannian metric on $H_n({\bf
R})$, because the sub-Riemannian distance is defined by the infimum
over a smaller class of curves.  However, it is also well-known that
the sub-Riemannian distance is topologically equivalent to the
Riemannian distance on $H_n({\bf R})$.

        Of course, this sub-Riemannian distance on $H_n({\bf R})$
depends only on the restriction of the given Riemannian metric to the
horizontal subspaces of the tangent spaces of $H_n({\bf R})$.  If the
Riemannian metric on $H_n({\bf R})$ is invariant under translations on
the right, as before, then the associated Riemannian and
sub-Riemannian distances on $H_n({\bf R})$ will have the same
property.  In the sub-Riemannian case, this uses the fact that the
class of horizontal curves in $H_n({\bf R})$ is invairiant under
translations on the right, because of the analogous property of
horizontal subspaces of the tangent spaces of $H_n({\bf R})$.

        An advantage of sub-Riemannian structures on $H_n({\bf R})$ is
that they behave more simply with respect to the dilations $\delta_r$.
Observe that a dilation $\delta_r$ applied to a horizontal curve in
$H_n({\bf R})$ is also a horizontal curve, because the horizontal
subspaces of the tangent spaces of $H_n({\bf R})$ are invariant under
$\delta_r$.  If we use a Riemannian metric on $H_n({\bf R})$ that is
invariant under right translations, then the restriction of this
Riemannian metric to the horizontal subspaces of the tangent spaces of
$H_n({\bf R})$ transforms by scalar multiplication under the dilations
$\delta_r$.  This implies that sub-Riemannian distances on $H_n({\bf
R})$ also transform by scalar multiplication under a dilation
$\delta_r$ under these conditions.  More precisely, the sub-Riemannian
distance between the images of two points in $H_n({\bf R})$ under
$\delta_r$ is equal to the absolute value $|r|$ of $r \in {\bf R}$
times the sub-Riemannian distance between the two points.

        Let us continue to use a Riemannian metric on $H_n({\bf R})$
that is invariant under right translations.  It is well known that the
corresponding sub-Riemannian distance between a point $(x, y, t) \in
H_n({\bf R})$ and $(0, 0, 0)$ is comparable to
\begin{equation}
\label{|x| + |y| + |t|^{1/2}}
        |x| + |y| + |t|^{1/2},
\end{equation}
where $|x|$, $|y|$ are the standard Euclidean norms of $x, y \in {\bf
R}^n$.  This means that the sub-Riemannian distance and (\ref{|x| +
|y| + |t|^{1/2}}) are each bounded by constant multiples of the other,
where the constants do not depend on the point $(x, y, t)$.  One can
estimate the sub-Riemannian distance between any two elements of
$H_n({\bf R})$ using invariance under right translations to reduce to
the case where one of the points is $(0, 0, 0)$.

        Because the distribution of horizontal subspaces of the
tangent spaces of $H_n({\bf R})$ is invariant under right
translations, there is an analogous distribution of horizontal
subspaces of the tangent spaces of $H_n({\bf R}) / H_n({\bf Z})$.
More precisely, the differential of the quotient mapping from
$H_n({\bf R})$ onto $H_n({\bf R}) / H_n({\bf Z})$ sends horizontal
subspaces of tangent spaces of the domain to horizontal subspaces of
the corresponding tangent spaces of the range.  Similar remarks apply
to other quotients of $H_n({\bf R})$ by discrete subgroups on the
right, and one can also consider the correspoinding sub-Riemannian
structures on these quotients.

\section{Invariant measures}
\label{invariant measures}
\setcounter{equation}{0}

        As usual, we can identify $H_n({\bf R})$ with ${\bf R}^n
\times {\bf R}^n \times {\bf R}$ as a smooth manifold, so that
$H_n({\bf R})$ is equipped with $(2n + 1)$-dimensional Lebesgue
measure in particular.  It is well known and not difficult to check
that Lebesgue measure on $H_n({\bf R})$ is invariant under
translations on the right and on the left with respect to $\diamond$,
using the standard formula for computing volumes under a change of
variables.  The main point is that the matrices corresponding to the
differentials of these translation mappings are upper-triangular with
$1$'s along the diagonal, so that their determinants are equal to $1$.
Equivalently, the standard volume form on ${\bf R}^n \times {\bf R}^n
\times {\bf R}$ of degree $2n + 1$ is invariant under right and left
translations with respect to $\diamond$.

        Using invariance of the volume form on $H_n({\bf R})$ under
right translations, we get a volume form on $H_n({\bf R}) / H_n({\bf
Z})$ too.  This volume form on $H_n({\bf R}) / H_n({\bf Z})$ is
characterized by the property that its pull-back under the quotient
mapping from $H_n({\bf R})$ onto $H_n({\bf R}) / H_n({\bf Z})$ is
equal to the original volume form on $H_n({\bf R})$.  Using invariance
of the volume form on $H_n({\bf R})$ under left translations, we get
that the volume form on $H_n({\bf R}) / H_n({\bf Z})$ is invariant
under under the natural action of $H_n({\bf R})$ on the quotient by
multiplication on the left.  Similar remarks apply to other quotients
of $H_n({\bf R})$ by discrete subgroups.

        Suppose that $H_n({\bf R})$ is equipped with a Riemannian
metric invariant under right translations, and consider the
corresponding sub-Riemannian distance on $H_n({\bf R})$, as in the
previous section.  If $B_r$ is a ball in $H_n({\bf R})$ with respect
to the sub-Riemannian distance with radius $r$, then the Lebesgue
measure of $B_r$ is equal to a positive constant times $r^{2 n + 2}$.
To see this, one can first use invariance of the sub-Riemannian
distance and of the volume under translations on the right to reduce
to the case where $B_r$ is centered at $(0, 0, 0)$.  One can then use
the dilation $\delta_r$ to reduce to the case where $r = 1$, since it
is easy to compute the behavior of Lebesgue measure under this
dilation.

        Let $k$ be a positive integer, and let $H_n({\bf Z} / k \,
{\bf Z})$ be the group defined as follows.  As a set, $H_n({\bf Z} / k
\, {\bf Z})$ is the same as
\begin{equation}
\label{(Z / k Z)^n times (Z / k Z)^n times (Z / k Z)}
        ({\bf Z} / k \, {\bf Z})^n \times ({\bf Z} / k \, {\bf Z})^n
                                        \times ({\bf Z} / k \, {\bf Z}).
\end{equation}
The group structure on $H_n({\bf Z} / k \, {\bf Z})$ can be defined in
the same way as for $H_n({\bf R})$, using the standard ring structure
on ${\bf Z} / k \, {\bf Z}$ as a quotient of ${\bf Z}$.  The usual
quotient mapping from ${\bf Z}$ onto ${\bf Z} / k \, {\bf Z}$ leads to
a homomorphism from $H_n({\bf Z})$ onto $H_n({\bf Z} / k \, {\bf Z})$
with kernel equal to $H_n(k \, {\bf Z})$.  In particular, the index of
$H_n(k \, {\bf Z})$ as a subgroup of $H_n({\bf Z})$ is equal to the
number of elements of $H_n({\bf Z} / k \, {\bf Z})$, which is $k^{2n +
1}$.  Similarly, if $r \ge 2$ is an integer, then $H_n(r \, {\bf Z}) /
\delta_r(H_n({\bf Z}))$ is isomorphic as a group to $r \, {\bf Z} /
r^2 \, {\bf Z} \cong {\bf Z} / r \, {\bf Z}$, which has $r$ elements.
It follows that the index of $\delta_r(H_n({\bf Z}))$ as a subgroup of
$H_n({\bf Z})$ is equal to the product of the index of
$\delta_r(H_n({\bf Z}))$ in $H_n(r \, {\bf Z})$ with the index of
$H_n(r \, {\bf Z})$ in $H_n({\bf Z})$, which is $r^{2 n + 2}$.

\section{$r$-Adic integers}
\label{r-adic integers}
\setcounter{equation}{0}

        Let $r \ge 2$ be an integer, and let $|a|_r$ be the $r$-adic
absolute value of an integer $a$, defined as follows.  If $a = 0$,
then we put $|a|_r = 0$.  Otherwise, if $a \ne 0$, then we let $l$ be
the largest nonnegative integer such that $a$ is an integer multiple
of $r^l$, and we put
\begin{equation}
\label{|a|_r = r^{-l}}
        |a|_r = r^{-l}.
\end{equation}
If $r = p$ is prime, then this is the same as the usual $p$-adic
absolute value of $a$.  It is easy to check that
\begin{equation}
\label{|a + b|_r le max(|a|_r, |b|_r)}
        |a + b|_r \le \max(|a|_r, |b|_r)
\end{equation}
and
\begin{equation}
\label{|a b|_r le |a|_r |b|_r}
        |a \, b|_r \le |a|_r \, |b|_r
\end{equation}
for every $a, b \in {\bf Z}$, even when $r$ is not prime.  If $r = p$
is prime, then equality holds in (\ref{|a b|_r le |a|_r |b|_r}).  In
this case, the $p$-adic absolute value can be defined on the field
${\bf Q}$ of rational numbers, with the same properties.

        The $r$-adic distance $d_r(a, b)$ is defined on ${\bf Z}$ by
\begin{equation}
\label{d_r(a, b) = |a - b|_r}
        d_r(a, b) = |a - b|_r.
\end{equation}
Observe that
\begin{equation}
\label{d_r(a, c) le max(d_r(a, b), d_r(b, c))}
        d_r(a, c) \le \max(d_r(a, b), d_r(b, c))
\end{equation}
for every $a, b, c \in {\bf Z}$, because of (\ref{|a + b|_r le
max(|a|_r, |b|_r)}).  Thus $d_r(a, b)$ defines an ultrametric on ${\bf
Z}$, which means that $d_r(a, b)$ defines a metric on ${\bf Z}$ that
satisfies the stronger ultrametric version (\ref{d_r(a, c) le
max(d_r(a, b), d_r(b, c))}) version of the triangle inequality.  If $r
= p$ is prime, then one gets an ultrametric on ${\bf Q}$ in the same
way.

        The $r$-adic integers ${\bf Z}_r$ may be defined as the
completion of ${\bf Z}$ as a metric space with respect to the $r$-adic
metric.  Addition and multiplication on ${\bf Z}$ can be extended to
${\bf Z}_r$, so that ${\bf Z}_r$ becomes a commutative ring.  The
$r$-adic absolute value and metric can also be extended from ${\bf Z}$
to ${\bf Z}_r$, with the same proprties as before.  If $r = p$ is
prime, then the completion of ${\bf Q}$ with respect to the $p$-adic
metric leads to the field ${\bf Q}_p$ of $p$-adic numbers.  In this
case, one can show that the $p$-adic integers are the same as the
$p$-adic numbers with $p$-adic absolute value less than or equal to
$1$.

        The set ${\bf Z}_r$ of $r$-adic integers is compact with
respect to the $r$-adic metric.  To see this, it suffices to show that
${\bf Z}_r$ is totally bounded with respect to the $r$-adic metric,
which means that ${\bf Z}_r$ can be covered by finitely many balls
with arbitarily small radius.  This uses the well-known fact that a
subset of a complete metric space is compact if it is closed and
totally bounded.  In this case, ${\bf Z}_r$ can be expressed as the
union of $r^l$ closed balls of radius $r^{-l}$ for each nonnegative
integer $l$.  More precisely, it suffices to use closed balls in ${\bf
Z}_r$ of radius $r^{-1}$ centered at representatives of the cosets of
$r^l \, {\bf Z}$ in ${\bf Z}$.

\section{Cartesian products}
\label{cartesian products}
\setcounter{equation}{0}

        Let $r \ge 2$ be an integer again, and consider the Cartesian product
\begin{equation}
\label{X = prod_{l = 1}^infty ({bf Z} / r^l {bf Z})}
        X = \prod_{l = 1}^\infty ({\bf Z} / r^l \, {\bf Z}).
\end{equation}
Thus $X$ consists of the sequences $x = \{x_l\}_{l = 1}^\infty$, where
$x_l$ is an element of the quotient ${\bf Z} / r^l \, {\bf Z}$ of
integers modulo $r^l$ for each $l \ge 1$.  Of course, ${\bf Z} / r^l
\, {\bf Z}$ is a commutative ring for each $l \ge 1$, with respect to
addition and multiplication of integers modulo $r^l$.  It follows that
$X$ is also a commutative ring with respect to coordinatewise addition
and multiplication.  We can consider $X$ as a Hausdorff topological
space as well, with respect to the product topology associated to the
discrete topology on ${\bf Z} / r^l \, {\bf Z}$ for each $l$.  Note
that $X$ is compact with respect to this topology, because ${\bf Z} /
r^l \, {\bf Z}$ is a finite set for each $l$.  It is easy to see that
addition and multiplication on $X$ define continuous mappings from $X
\times X$ into $X$ with respect to this topology, using the product
topology on $X \times X$ corresponding to the topology just defined on
$X$.  Similarly, the additive inverse $-x$ of $x \in X$ corresponds to
the sequence $\{-x_l\}_{l = 1}^\infty$ in $X$, and $x \mapsto -x$
is a continuous mapping from $X$ into itself.  This shows that $X$
is a topological ring with respect to this topology.

        If $x$, $y$ are distinct elements of $X$, then let $l(x, y)$
be the smallest positive integer $l$ such that $x_l \ne y_l$, which is
the same as the largest positive integer such that $x_l = y_l$ for
every $l < l(x, y)$.  Put
\begin{equation}
\label{rho(x, y) = r^{-l(x, y) + 1}}
        \rho(x, y) = r^{-l(x, y) + 1}
\end{equation}
in this case, and $\rho(x, y) = 0$ when $x = y$, which amounts to
taking $l(x, y) = +\infty$ when $x = y$.  It is easy to see that
\begin{equation}
\label{l(x, z) ge min(l(x, y), l(y, z))}
        l(x, z) \ge \min(l(x, y), l(y, z))
\end{equation}
and hence that
\begin{equation}
\label{rho(x, z) le max(rho(x, y), rho(y, z))}
        \rho(x, z) \le \max(\rho(x, y), \rho(y, z))
\end{equation}
for every $x, y, z \in X$.  This implies that $\rho(x, y)$ defines an
ultrametric on $X$, which is to say that $\rho(x, y)$ is a metric on
$X$ that satisfies the stronger ultrametric version (\ref{rho(x, z) le
max(rho(x, y), rho(y, z))}) of the triangle inequality.  Note that the
topology on $X$ determined by this ultrametric is the same as the one
described in the previous paragraph.

        A sequence $x(1), x(2), x(3), \ldots$ of elements of $X$
converges to an element $x$ of $X$ if and only if for each positive
integer $l$, $x_l(j) = x_l$ for all sufficiently large $j$, depending
on $l$.  Similarly, $x(1), x(2), x(3), \ldots$ is a Cauchy sequence in
$X$ with respect to $\rho(\cdot, \cdot)$ if and only if for each
positive integer $l$, $x_l(j)$ is eventually constant in $j$,
depending on $l$.  Using this, it is easy to see that every Cauchy
sequence in $X$ with respect to $\rho(\cdot, \cdot)$ converges to an
element of $X$, so that $X$ is complete as a metric space with respect
to $\rho(\cdot, \cdot)$.  Of course, compact metric spaces are always
complete, but in this case we can check this more directly.

\section{Coherent sequences}
\label{coherent sequences}
\setcounter{equation}{0}

        Let us continue with the same notation and hypotheses as in
the preceding section.  Because $r^{l + 1} \, {\bf Z} \subseteq r^l \,
{\bf Z}$, there is a natural ring homomorphism from ${\bf Z} / r^{l +
1} \, {\bf Z}$ onto ${\bf Z} / r^l \, {\bf Z}$ for each $l \ge 1$.  An
element $x = \{x_l\}_{l = 1}^\infty$ of $X$ is said to be a
\emph{coherent sequence} if $x_l$ is the image in ${\bf Z} / r^l \,
{\bf Z}$ of $x_{l + 1} \in {\bf Z} / r^{l + 1} \, {\bf Z}$ for each
$l$.  Let $Y$ be the subset of $X$ consisting of coherent sequences.
It is easy to see that $Y$ is a subring of $X$ with respect to
coordinatewise addition and multiplication, and that $Y$ is a closed
subset of $X$ with respect to the product topology on $X$.

        Let $q_l$ be the natural quotient mapping from ${\bf Z}$ onto
${\bf Z} / r^l \, {\bf Z}$ for each $l$, and let $q$ be the mapping
from ${\bf Z}$ into $X$ defined by
\begin{equation}
        q(a) = \{q_l(a)\}_{l = 1}^\infty
\end{equation}
for each $a \in {\bf Z}$.  Thus $q$ is a ring homomorphism, since
$q_l$ is a ring homomorphism for each $l$.  Note that the kernel of
$q$ in ${\bf Z}$ is trivial, so that $q$ is injective.  It is easy to
see that $q(a)$ is a coherent sequence for each $a \in {\bf Z}$,
because $q_l$ is the same as the composition of $q_{l + 1}$ with the
natural homomorphism from ${\bf Z} / r^{l + 1} \, {\bf Z}$ onto ${\bf
Z} / r^l \, {\bf Z}$ for each $l$.  One can also check that $q({\bf
Z})$ is dense in $Y$, so that $Y$ is the same as the closure of
$q({\bf Z})$ in the product topology on $X$.  More precisely, if $x =
\{x_l\}_{l = 1}^\infty$ is any coherent sequence and $L$ is a positive
integer, then one can choose $a \in {\bf Z}$ such that $q_L(a) = x_L$.
This implies that
\begin{equation}
        q_l(a) = x_l
\end{equation}
for every $l \le L$, because $q(a)$ and $x$ are coherent sequences,
which says exactly that $x$ can be approximated by elements $q(a)$ of
$q({\bf Z})$ with respect to the product topology on $X$, as desired.

        If $\rho(x, y)$ is the ultrametric on $X$ defined in the
previous section, then
\begin{equation}
\label{rho(q(a), q(b)) = d_r(a, b) = |a - b|_r}
        \rho(q(a), q(b)) = d_r(a, b) = |a - b|_r
\end{equation}
for every $a, b \in {\bf Z}$, where $|a|_r$ and $d_r(a, b)$ are the
$r$-adic absolute value and $r$-adic metric, as in Section \ref{r-adic
integers}.  This is trivial when $a = b$, and otherwise $l(q(a),
q(b))$ is defined in Section \ref{cartesian products} as the smallest
positive integer $l$ such that $q_l(a) \ne q_l(b)$.  This is the same
as the smallest positive integer $l$ such that $q_l(a - b) \ne 0$, so
that $l_1 = l(q(a), q(b)) - 1$ is the largest nonnegative integer such
that $a - b$ is an integer multiple of $r^{l_1}$.  Thus
(\ref{rho(q(a), q(b)) = d_r(a, b) = |a - b|_r}) follows from
(\ref{|a|_r = r^{-l}}) and (\ref{rho(x, y) = r^{-l(x, y) + 1}}).

        This shows that $q$ defines an isometric embedding of ${\bf
Z}$ with the $r$-adic metric into $X$ with the ultrametric $\rho(x,
y)$.  Because $X$ is complete as a metric space with respect to
$\rho(x, y)$, the completion ${\bf Z}_r$ of ${\bf Z}$ with respect to
the $r$-adic metric can be identified with the closure of $q({\bf Z})$
in $X$, which we have seen is the same as the set $Y$ of coherent
sequences in $X$.  It is easy to see that this identification is also
compatible with addition and multiplication of $r$-adic integers and
coherent sequences.  Equivalently, if $\{a_j\}_{j = 1}^\infty$ is a
Cauchy sequence of integers with respect to the $r$-adic metric, then
for each positive integer $l$, there is an $x_l \in {\bf Z} / r^l \,
{\bf Z}$ such that $q_l(a_j) = x_l$ for all sufficiently large $j$,
depending on $l$.  One can check that $x = \{x_l\}_{l = 1}^\infty$ is
a coherent sequence under these conditions, and that this corresponds
to the identification between ${\bf Z}_r$ and $Y$ just mentioned.

\section{A related construction}
\label{related construction}
\setcounter{equation}{0}

        Let $r \ge 2$ be an integer again, and consider the Cartesian product
\begin{equation}
\label{widetilde{X} = prod_{l = 0}^infty ({bf R} / r^l {bf Z})}
        \widetilde{X} = \prod_{l = 0}^\infty ({\bf R} / r^l \, {\bf Z}).
\end{equation}
Here ${\bf R} / r^l \, {\bf Z}$ refers to the quotient of ${\bf R}$ as
a commutative group with respect to addition by the subgroup $r^l \,
{\bf Z}$.  Thus $\widetilde{X}$ consists of the sequences $x =
\{x_l\}_{l = 0}^\infty$ with $x_l \in {\bf R} / r^l \, {\bf Z}$ for
each $l$, and $\widetilde{X}$ is also a commutative group with respect
to coordinatewise addition.  As before, ${\bf R} / r^l \, {\bf Z}$ is
a topological space in a natural way, such that the usual quotient
mapping $\widetilde{q}_l$ from ${\bf R}$ onto ${\bf R} / r^l \, {\bf
Z}$ is a local homeomorphism for each $l$.  Using the corresponding
product topology on $\widetilde{X}$, we get a compact Hausdorff space
which is homeomorphic to the product of circles.  It is easy to see
that $\widetilde{X}$ is a topological group, which is to say that the
group operations are continuous with respect to this topology.  The
space $X$ described in Section \ref{cartesian products} may be
considered as a subset of $\widetilde{X}$, by extending each sequence
$x = \{x_l\}_{l = 1}^\infty \in X$ to a sequence starting at $l = 0$,
where $x_0 = 0$ in ${\bf R} / {\bf Z}$.  With this identification, $X$
corresponds to a closed subgroup of $\widetilde{X}$.

        As in Section \ref{coherent sequences}, there is a natural
group homomorphism from ${\bf R} / r^{l + 1} \, {\bf Z}$ onto ${\bf R}
/ r^l \, {\bf Z}$ for each $l \ge 0$, because $r^l \, {\bf Z}
\subseteq r^{l + 1} \, {\bf Z}$.  An element $x = \{x_l\}_{l =
0}^\infty$ of $\widetilde{X}$ is said to be a \emph{coherent sequence}
if $x_l$ is the image in ${\bf R} / r^l \, {\bf Z}$ of $x_{l + 1} \in
{\bf R} / r^{l + 1} \, {\bf Z}$ for each $l \ge 0$.  Let
$\widetilde{Y}$ be the set of coherent sequences in $\widetilde{X}$,
and observe that $\widetilde{Y}$ is a closed subgroup of
$\widetilde{X}$.  Using the identification of $X$ with a closed
subgroup of $\widetilde{X}$ mentioned in the preceding paragraph, the
set $Y$ of coherent sequences in $X$ corresponds exactly to the
intersection of $\widetilde{Y}$ with $X$.  Remember that $Y$ can also
be identified with the set ${\bf Z}_r$ of $r$-adic integers, as in the
previous section.

        Let $\widetilde{q}$ be the mapping from ${\bf R}$ into
$\widetilde{X}$ defined by
\begin{equation}
\label{widetilde{q}(a) = {widetilde{q}_l(a)}_{l = 0}^infty}
        \widetilde{q}(a) = \{\widetilde{q}_l(a)\}_{l = 0}^\infty
\end{equation}
for each $a \in {\bf R}$, where $\widetilde{q}_l$ is the natural
quotient mapping from ${\bf R}$ onto ${\bf R} / r^l \, {\bf Z}$
mentioned earlier.  Note that $\widetilde{q}$ is a group homomorphism
from ${\bf R}$ into $\widetilde{X}$ with trivial kernel, and that
$\widetilde{q}$ is also a continuous mapping from ${\bf R}$ into
$\widetilde{X}$, because $\widetilde{q}_l$ is continuous for each $l$.
As in Section \ref{coherent sequences}, $q(a)$ is a coherent sequence
in $\widetilde{X}$ for each $a \in {\bf R}$, because $\widetilde{q}_l$
is the same as the composition of $\widetilde{q}_{l + 1}$ with the
natural homomorphism from ${\bf R} / r^{l + 1} \, {\bf Z}$ onto ${\bf
R} / r^l \, {\bf Z}$ for each $l$.  One can also check that
$\widetilde{q}({\bf R})$ is dense in $\widetilde{Y}$, so that
$\widetilde{Y}$ is the closure of $\widetilde{q}({\bf R})$ in
$\widetilde{X}$ with respect to the product topology, for essentially
the same reasons as before.

        Let $\pi_0$ be the coordinate projection from $\widetilde{X}$
onto ${\bf R} / {\bf Z}$, so that
\begin{equation}
\label{pi_0(x) = x_0}
        \pi_0(x) = x_0
\end{equation}
for each $x = \{x_l\}_{l = 0}^\infty \in \widetilde{X}$.  Thus $\pi_0$
is a continuous group homomorphism from $\widetilde{X}$ onto ${\bf R}
/ {\bf Z}$.  The restriction of $\pi_0$ to $\widetilde{Y}$ is a
continuous group homomorphism from $\widetilde{Y}$ onto ${\bf R} /
{\bf Z}$ with kernel equal to $Y$.

\section{$r$-Adic Heisenberg groups}
\label{r-adic heisenberg groups}
\setcounter{equation}{0}

        Let $r \ge 2$ be an integer again, and let ${\bf Z}_r$ be the
ring of $r$-adic integers, as in Section \ref{r-adic integers}.  The
corresponding Heisenberg group $H_n({\bf Z}_r)$ is defined first as a
set as ${\bf Z}_r^n \times {\bf Z}_r^n \times {\bf Z}_r$.  The group
structure $\diamond$ can be defined on $H_n({\bf Z}_r)$ exactly as in
(\ref{(x, y, t) diamond (x', y', t') = (x + x', y + y', t + t' + x
cdot y')}), and is associative for the same reasons as before.
Similarly, inverses in $H_n({\bf Z}_r)$ can be defined as in (\ref{(x,
y, t)^{-1} = (-x, -y, -t + x cdot y)}).  It is easy to see that the
group operations on $H_n({\bf Z}_r)$ are continuous with respect to
the topology on $H_n({\bf Z}_r)$ defined by the usual $r$-adic
topology on ${\bf Z}_r$ and the product topology, so that $H_n({\bf
Z}_r)$ is a topological group.  Note that $H_n({\bf Z}_r)$ is a
compact Hausdorff space with respect to this topology.  The Heisenberg
group $H_n({\bf Z})$ associated to the ordinary integers is a dense
subgroup of $H_n({\bf Z}_r)$, because ${\bf Z}$ is dense in ${\bf
Z}_r$.

        Consider the Cartesian product
\begin{equation}
\label{W = prod_{l = 1}^infty (H_n({bf Z}) / H_n(r^l {bf Z}))}
        W = \prod_{l = 1}^\infty (H_n({\bf Z}) / H_n(r^l \, {\bf Z})),
\end{equation}
so that the elements of $W$ are the sequences $w = \{w_l\}_{l = 1}^\infty$
such that $w_l$ is in $H_n({\bf Z}) / H_n(r^l \, {\bf Z})$ for each $l$.
Remember that $H_n(r^l \, {\bf Z})$ is a normal subgroup of $H_n({\bf Z})$
for each $l$, and that the quotient $H_n({\bf Z}) / H_n(r^l \, {\bf Z})$
is isomorphic to $H_n({\bf Z} / r^l \, {\bf Z})$.  Thus $W$ is also a
group, where the group operations are defined coordinatewise, and which
is isomorphic to
\begin{equation}
\label{prod_{l = 1}^infty H_n({bf Z} / r^l {bf Z})}
        \prod_{l = 1}^\infty H_n({\bf Z} / r^l \, {\bf Z}).
\end{equation}
If we equip $W$ with the product topology associated to the discrete
topology on $H_n({\bf Z}) / H_n(r^l \, {\bf Z})$ for each $l$, then
$W$ is a compact Hausdorff space, because $H_n({\bf Z}) / H_n(r^l \,
{\bf Z})$ has only finitely many elements for each $l$.  It is easy to
see that the group operations on $W$ are continuous with respect to
this topology, so that $W$ is a topological group.

        As usual, $H_n(r^{l + 1} \, {\bf Z}) \subseteq H_n(r^l \, {\bf
Z})$ for each $l$, which leads to a natural group homomorphism from
$H_n({\bf Z}) / H_n(r^{l + 1} \, {\bf Z})$ onto $H_n({\bf Z}) /
H_n(r^l \, {\bf Z})$ for each $l$.  An element $w = \{w_l\}_{l =
1}^\infty$ of $W$ is said to be a \emph{coherent sequence} if $w_l$ is
the image in $H_n({\bf Z}) / H_n(r^l \, {\bf Z})$ of $w_{l + 1} \in
H_n({\bf Z}) / H_n(r^{l + 1} \, {\bf Z})$ for each $l$.  Let $V$ be
the set of coherent sequences in $W$, which is a closed subgroup of
$W$.  Also let $\phi_l$ be the usual quotient mapping from $H_n({\bf
Z})$ onto $H_n({\bf Z}) / H_n(r^l \, {\bf Z})$, and let $\phi$ be the
mapping from $H_n({\bf Z})$ into $W$ whose $l$th component is equal to
$\phi_l$ for each $l$.  Thus $\phi$ is a group homomorphism from
$H_n({\bf Z})$ into $W$ with trivial kernel.  As before, $\phi$
actually maps $H_n({\bf Z})$ into $V$, because $\phi_l$ is the same as
the composition of $\phi_{l + 1}$ with the natural homomorphism from
$H_n({\bf Z}) / H_n(r^{l + 1} \, {\bf Z})$ onto $H_n({\bf Z}) /
H_n(r^l \, {\bf Z})$.  One can check that $\phi(H_n({\bf Z}))$ is
dense in $V$, so that $V$ is the same as the closure of $\phi(H_n({\bf
Z}))$ in $W$.

        The main point now is that $V$ is isomorphic to the Heisenberg
group $H_n({\bf Z}_r)$ associated to the $r$-adic numbers.  One way to
see this is to show that $\phi$ can be extended to a continuous
mapping from $H_n({\bf Z}_r)$ onto $V$, and that this extension is an
isomorphism.  Alternatively, one can use the identification of
$r$-adic integers with coherent sequences in (\ref{X = prod_{l =
1}^infty ({bf Z} / r^l {bf Z})}) discussed in Section \ref{coherent
sequences}.  In this case, it is helpful to remember the
identification of $W$ with (\ref{prod_{l = 1}^infty H_n({bf Z} / r^l
{bf Z})}), and to identify $V$ with the set of coherent sequences in
(\ref{prod_{l = 1}^infty H_n({bf Z} / r^l {bf Z})}).  This uses the
natural homomorphism from $H_n({\bf Z} / r^{l + 1} \, {\bf Z})$ onto
$H_n({\bf Z} / r^l \, {\bf Z})$ for each $l$, which comes from the
natural homomorphism from ${\bf Z} / r^{l + 1} \, {\bf Z}$ onto ${\bf
Z} / r^l \, {\bf Z}$.

\section{Heisenberg solenoids}
\label{heisenberg solenoids}
\setcounter{equation}{0}

        Let $r \ge 2$ be an integer again, and consider the Cartesian product
\begin{equation}
\label{widetilde{W} = prod_{l = 0}^infty (H_n({bf R}) / H_n(r^l {bf Z}))}
 \widetilde{W} = \prod_{l = 0}^\infty (H_n({\bf R}) / H_n(r^l \, {\bf Z})).
\end{equation}
As before, $H_n({\bf R}) / H_n(r^l \, {\bf Z})$ may be considered as a
compact smooth manifold of dimension $2 n + 1$ for each $l \ge 0$, so
that $\widetilde{W}$ is a compact Hausdorff space with respect to the
corresponding product topology.  There is also a natural action of
$H_n({\bf R})$ on $H_n({\bf R}) / H_n(r^l \, {\bf Z})$ by
multiplication on the left for each $l$, which leads to a continuous
action of $H_n({\bf R})$ on $\widetilde{W}$.  The space $W$ discussed
in the previous section may be identified with a subset of $\widetilde{W}$,
by extending each sequence $w = \{w_l\}_{l = 1}^\infty \in W$ to $l =
0$, with $w_0$ equal to the element of $H_n({\bf R}) / H_n({\bf Z})$
corresponding to $H_n({\bf Z}) / H_n({\bf Z})$.  Note that $W$ corresponds
to a closed set in $\widetilde{W}$, and that the topology on $W$ induced
from the one on $\widetilde{W}$ is the same as the topology on $W$ defined
in the previous section.

        As usual, because $H_n(r^{l + 1} \, {\bf Z}) \subseteq H_n(r^l
\, {\bf Z})$, there is a natural mapping from $H_n({\bf R}) / H_n(r^{l
+ 1} \, {\bf Z})$ onto $H_n({\bf R}) / H_n(r^l \, {\bf Z})$ for each
$l \ge 0$, and this mapping is a local diffeomorphism.  An element $w
= \{w_l\}_{l = 0}^\infty$ of $\widetilde{W}$ is said to be a
\emph{coherent sequence} if $w_l$ is the image in $H_n({\bf R}) /
H_n(r^l \, {\bf Z})$ of $w_{l + 1} \in H_n({\bf R}) / H_n(r^{l + 1} \,
{\bf Z})$ for each $l$.  Let $\widetilde{V}$ be the set of coherent
sequences in $\widetilde{W}$, which is a closed set in $\widetilde{W}$
with respect to the product topology that is invariant under the
action of $H_n({\bf R})$ on $\widetilde{W}$ by multiplication on the
left mentioned earlier.  Using the identification of $W$ with a closed
set in $\widetilde{W}$ described in the preceding paragraph, the set
$V$ of coherent sequences in $W$ corresponds to the intersection of
$\widetilde{V}$ and $W$.  Remember that $V$ can also be identified
with $H_n({\bf Z}_r)$, as in the previous section.

        Let $\widetilde{\phi}_l$ be the natural quotient mapping from
$H_n({\bf R})$ onto $H_n({\bf R}) / H_n(r^l \, {\bf Z})$ for each $l
\ge 0$, which is a local diffeomorphism.  This leads to a continuous
embedding $\widetilde{\phi}$ of $H_n({\bf R})$ into $\widetilde{W}$,
whose $l$th component is equal to $\widetilde{\phi}_l$ for each $l$.
Note that this embedding intertwines the natural actions of $H_n({\bf
R})$ on itself and on $\widetilde{W}$ by multiplication on the left.
As usual, $\widetilde{\phi}$ actually maps $H_n({\bf R})$ into
$\widetilde{V}$, because $\widetilde{\phi}_l$ is the same as the
composition of $\widetilde{\phi}_{l + 1}$ with the natural mapping
from $H_n({\bf R}) / H_n(r^{l + 1} \, {\bf Z})$ onto $H_n({\bf R}) /
H_n(r^l \, {\bf Z})$ for each $l$.  One can also check that
$\widetilde{\phi}(H_n({\bf R}))$ is dense in $\widetilde{V}$ with
respect to the product topology on $\widetilde{W}$, so that
$\widetilde{V}$ is the closure of $\widetilde{\phi}(H_n({\bf R}))$ in
$\widetilde{W}$.

\section{Another version}
\label{another version}
\setcounter{equation}{0}

        Let us continue with the same notation and hypotheses as in
the preceding section, and consider the Cartesian product
\begin{equation}
\label{widetilde{U} = prod_{l = 0}^infty (H_n(R) / delta_{r^l}(H_n(Z)))}
        \widetilde{U} = \prod_{l = 0}^\infty
                          (H_n({\bf R}) / \delta_{r^l}(H_n({\bf Z}))).
\end{equation}
As usual, $H_n({\bf R}) / \delta_{r^l}(H_n({\bf Z}))$ may be
considered as a compact smooth manifold of dimension $2 n + 1$ for
each $l$, so that $\widetilde{U}$ is a compact Hausdorff space with
respect to the product topology.  There is also a natural continuous
action of $H_n({\bf R})$ on $\widetilde{U}$, using the natural actions
of $H_n({\bf R})$ on $H_n({\bf R}) / \delta_{r^l}(H_n({\bf Z}))$ by
multiplication on the left for each $l$.

        As before, because $\delta_{r^{l + 1}}(H_n({\bf Z})) \subseteq
\delta_{r^l}(H_n({\bf Z}))$, there is a natural mapping from $H_n({\bf
R}) / \delta_{r^{l + 1}}(H_n({\bf Z}))$ onto $H_n({\bf R}) /
\delta_{r^l}(H_n({\bf Z}))$ for each $l \ge 0$, which is a local
diffeomorphism.  An element $u = \{u_l\}_{l = 0}^\infty$ of
$\widetilde{U}$ is said to be a \emph{coherent sequence} if $u_l$ is
the image in $H_n({\bf R}) / \delta_{r^l}(H_n({\bf Z}))$ of $u_{l + 1}
\in H_n({\bf R}) / \delta_{r^{l + 1}}(H_n({\bf Z}))$ for each $l$.
The set $\widetilde{V}'$ of coherent sequences in $\widetilde{U}$ is a
closed set in $\widetilde{U}$ with respect to the product topology
that is invariant under the action of $H_n({\bf R})$ on
$\widetilde{U}$ by multiplication on the left mentioned earlier.  
Let $\widetilde{\psi}_l$ be the natural quotient mapping from
$H_n({\bf R})$ onto $H_n({\bf R}) / \delta_{r^l}(H_n({\bf Z}))$ for
each $l \ge 0$, which is a local diffeomorphism that intertwines the
actions of $H_n({\bf R})$ on itself and on $H_n({\bf R}) /
\delta_{r^l}(H_n({\bf Z}))$ by multiplication on the left.
This leads to a continuous embedding $\widetilde{\psi}$ of $H_n({\bf
R})$ into $\widetilde{U}$, whose $l$th component is equal to
$\widetilde{\psi}_l$ for each $l$, and which intertwines the actions
of $H_n({\bf R})$ on itself and on $\widetilde{U}$ by multiplication
on the left.  This embedding actually maps $H_n({\bf R})$ into the set
$\widetilde{V}'$ of coherent sequences, because $\widetilde{\psi}_l$
is the same as the composition of $\widetilde{\psi}_{l + 1}$ with the
natural mapping from $H_n({\bf R}) / \delta_{r^{l + 1}}(H_n({\bf Z}))$
onto $H_n({\bf R}) / \delta_{r^l}(H_n({\bf Z}))$ for each $l$.
One can check that $\widetilde{\psi}(H_n({\bf R}))$ is dense in
$\widetilde{V}'$ with respect to the product topology on $\widetilde{U}$,
so that $\widetilde{V}'$ is the same as the closure of 
$\widetilde{\psi}(H_n({\bf R}))$ in $\widetilde{U}$.

        There is a natural identification between the set
$\widetilde{V}'$ of coherent sequences in $\widetilde{U}$ and the set
$\widetilde{V}$ of coherent sequences in $\widetilde{W}$, because of
the way that the subgroups $H_n(r^l \, {\bf Z})$ and
$\delta_{r^l}(H_n({\bf Z}))$ of $H_n({\bf Z})$ are interlaced, as in
(\ref{H_n(r^2 {bf Z}) subseteq delta_r(H_n({bf Z})) subseteq H_n(r {bf
Z})}).  With respect to this identification, the embeddings
$\widetilde{\phi}$ and $\widetilde{\psi}$ of $H_n({\bf R})$ into
$\widetilde{V}$ and $\widetilde{V}'$ are the same, the actions of
$H_n({\bf R})$ on $\widetilde{V}$ and $\widetilde{V}'$ by
multiplication on the left are the same, the topologies on
$\widetilde{V}$ and $\widetilde{V}'$ are the same, and the projections
of $\widetilde{V}$ and $\widetilde{V}'$ onto $H_n({\bf R}) / H_n({\bf
Z})$ from the $l = 0$ coordinates are the same.

        Remember that $H_n({\bf R}) / \delta_{r^l}(H_n({\bf Z}))$ is
diffeomorphic to $H_n({\bf R}) / H_n({\bf Z})$ for each $l$, because
$\delta_{r^l}$ is an automorphism of $H_n({\bf R})$.  This permits us
to identify $\widetilde{U}$ with the Cartesian product of a sequence
of copies of $H_n({\bf R}) / H_n({\bf Z})$.  This identification
affects the way that $H_n({\bf R})$ acts by multiplication on the
left, but this is easy to track.  Using these diffeomorphisms, the
natural mappings from $H_n({\bf R}) / \delta_{r^{l + 1}}(H_n({\bf
Z}))$ onto $H_n({\bf R}) / \delta_{r^l}(H_n({\bf Z}))$ correspond to
the same local diffeomorphism from $H_n({\bf R}) / H_n({\bf Z})$ onto
itself.  This local diffeomorphism is defined by applying $\delta_r$
on $H_n({\bf R})$ and then passing to the quotient by $H_n({\bf Z})$
on the domain and range, which is possible because $\delta_r$ is an
automorphism of $H_n({\bf R})$ that maps $H_n({\bf Z})$ into itself.


\begin{thebibliography}{49}

\addcontentsline{toc}{section}{References}





\bibitem {a-m} M.~Atiyah and I.~MacDonald, {\it Introduction to
Commutative Algebra}, Addison-Wesley, 1969.

\bibitem {b} R.~Beals, {\it Analysis: An Introduction}, Cambridge
University Press, 2004.

\bibitem {b-g} R.~Beals and P.~Greiner, {\it Calculus on Heisenberg
Manifolds}, Princeton University Press, 1988.

\bibitem {b-r} A.~Bella\"iche and J.-J.~Risler, editors, {\it
Sub-Riemannian Geometry}, Birkh\"auser, 1996.

\bibitem {b-mac} G.~Birkhoff and S.~Mac Lane, {\it A Survey of Modern
Algebra}, 4th edition, Macmillan, 1977.

\bibitem {c-d-p-t} L.~Capogna, D.~Danielli, S.~Pauls, and J.~Tyson,
{\it An Introduction to the Heisenberg Group and the Sub-Riemannian
Isoperimetric Problem}, Birkh\"auser, 2007.

\bibitem {car} H.~Cartan, {\it Elementary Theory of Analytic Functions
of One or Several Complex Variables}, translated from the French,
Dover, 1995.

\bibitem {cas} J.~Cassels, {\it Local Fields}, Cambridge University
Press, 1986.

\bibitem {c-w-1} R.~Coifman and G.~Weiss, {\it Analyse Harmonique
Non-Commutative sur Certains Espaces Homog\`enes}, Lecture Notes in
Mathematics {\bf 242}, Springer-Verlag, 1971.

\bibitem {c-w-2} R.~Coifman and G.~Weiss, {\it Extensions of Hardy
spaces and their use in analysis}, Bulletin of the American
Mathematical Society {\bf 83} (1977), 569--645.

\bibitem {d-s-1} G.~David and S.~Semmes, {\it Fractured Fractals and
Broken Dreams: Self-Similar Geometry through Metric and Measure},
Oxford University Press, 1997.

\bibitem {gf} G.~Folland, {\it Real Analysis: Modern Techniques and
their Applications}, 2nd edition, Wiley, 1999.

\bibitem {f-s} G.~Folland and E.~Stein, {\it Hardy Spaces on
Homogeneous Groups}, Princeton University Press, 1982.

\bibitem {fg} F.~Gouv\^ea, {\it $p$-Adic Numbers: An Introduction},
2nd edition, Springer-Verlag, 1997.

\bibitem {h1} J.~Heinonen, {\it Calculus on Carnot groups}, in {\it
Fall School in Analysis (Jyv\"askyl\"a, 1994)}, 1--31, Reports of the
Department of Mathematics and Statistics {\bf 68}, University of
Jyv\"askyl\"a, 1995.

\bibitem {h2} J.~Heinonen, {\it Lectures on Analysis on Metric
Spaces}, Springer-Verlag, 2001.

\bibitem {jh} J.~Humphreys, {\it Introduction to Lie Algebras and
Representation Theory}, Springer-Verlag, 1978.

\bibitem {nj} N.~Jacobson, {\it Lie Algebras}, Dover, 1979.

\bibitem {k-r} A.~Koranyi and H.~Reimann, {\it Quasiconformal mappings
on the Heisenberg group}, Inventiones Mathematicae {\bf 80} (1985),
309--338.

\bibitem {sk1} S.~Krantz, {\it A Panorama of Harmonic Analysis},
Mathematical Association of America, 1999.

\bibitem {sk2} S.~Krantz, {\it Function Theory of Several Complex
Variables}, AMS Chelsea, 2001.

\bibitem {sk3} S.~Krantz, {\it Explorations in Harmonic Analysis},
with the assistance of L.~Lee, Birkh\"auser, 2009.

\bibitem {mac-b} S.~Mac Lane and G.~Birkhoff, {\it Algebra}, 3rd
edition, Chelsea, 1988.

\bibitem {mc} G.~McCarty, {\it Topology: An Introduction with
Application to Topological Groups}, 2nd edition, Dover, 1988.

\bibitem {rm} R.~Montgomery, {\it A Tour of Subriemannian Geometries,
their Geodesics and Applications}, American Mathematical Society,
2002.

\bibitem {r1} W.~Rudin, {\it Principles of Mathematical Analysis}, 3rd
edition, McGraw-Hill, 1976.

\bibitem {r2} W.~Rudin, {\it Real and Complex Analysis}, 3rd edition,
McGraw-Hill, 1987.

\bibitem {r3} W.~Rudin, {\it Fourier Analysis on Groups}, Wiley, 1990.

\bibitem {r4} W.~Rudin, {\it Function Theory on the Unit Ball in ${\bf
C}^n$}, Springer-Verlag, 2008.

\bibitem {s1} S.~Semmes, {\it An introduction to analysis on metric
spaces}, Notices of the American Mathematical Society {\bf 50} (2003),
438--443.

\bibitem {s2} S.~Semmes, {\it An introduction to Heisenberg groups
in analysis and geometry}, Notices of the American Mathematical
Society {\bf 50} (2003), 640--646.

\bibitem {jps-1} J.-P.~Serre, {\it Local Fields}, translated from the
French by M.~Greenberg, Springer-Verlag, 1979.

\bibitem {jps-2} J.-P.~Serre, {\it Lie Algebras and Lie Groups}, 2nd
edition, Lecture Notes in Mathematics {\bf 1500}, Springer-Verlag,
2006.

\bibitem {st} E.~Stein, {\it Harmonic Analysis: Real-Variable Methods,
Orthogonality, and Oscillatory Integrals}, with the assistance of
T.~Murphy, Princeton University Press, 1993.

\bibitem {s-w} E.~Stein and G.~Weiss, {\it Introduction to Fourier
Analysis on Euclidean Spaces}, Princeton University Press, 1971.

\bibitem {ds} D.~Sullivan, {\it Linking the universalities of
Milnor--Thurston, Feigenbaum and Ahlfors--Bers}, in {\it Topological
Methods in Modern Mathematics}, 543--564, Publish or Perish, 1993.

\bibitem {t} M.~Taibleson, {\it Fourier Analysis on Local Fields},
Princeton University Press, 1975.

\bibitem {th} W.~Thurston, {\it Three-Dimensional Geometry and
Topology}, edited by S.~Levy, Princeton University Press, 1997.

\bibitem {v-s-c} N.~Varopoulos, L.~Saloff-Coste, and T.~Coulhon, {\it
Analysis and Geometry on Groups}, Cambridge University Press, 1992.

\bibitem {w} A.~Weil, {\it Basic Number Theory}, Springer-Verlag,
1995.



\end{thebibliography}
\end{document}